\documentclass{article}

\usepackage[english]{babel}

\usepackage[letterpaper,top=2cm,bottom=2cm,left=3cm,right=3cm,marginparwidth=1.75cm]{geometry}

\usepackage{amsmath}
\usepackage{graphicx}
\usepackage[colorlinks=true, allcolors=blue]{hyperref}
\usepackage{booktabs}

\begin{document}

\begin{center}
{\bf Coupled Van der Pol Networks}
\end{center}

\begin{center}
Pokutnyi O. \footnote{Institute of mathematics of NAS of Ukraine, 01024, Tereshenkivska 3, Kiev, Ukraine, alex$\_$poker@imath.kive.ua}, Viatchaninov O. \footnote{Independent researcher, alex.viatchaninov@gmail.com}
\end{center}

$~~~$In this work, we provide a complete analytical characterization of amplitude modes in a nonlinear Van der Pol network. 

The phase condition $c_{1m}^{k} = \pm c_{2m}^{k}$ fully classifies all possible configurations, leading to a finite set of discrete solutions.

This result provides a complete description of the system’s limit cycles, including all stable oscillatory states for arbitrary network size, without relying on approximate or iterative numerical schemes.

We also derive necessary and sufficient conditions for the existence of solutions of a coupled system of Van der Pol equations in a Hilbert space with periodic boundary conditions. These conditions depend on the oscillation amplitudes of the associated generating problem.

Such systems form a paradigmatic class exhibiting self-sustained oscillations. We extend classical amplitude analysis to infinite-dimensional coupled systems and uncover invariant structures that remain robust under changes in network size. We also highlight connections with quantum geometric phases~\cite{Berry1984, Simon1983}.

In the context of the van der Pol oscillator, the obtained result can be interpreted as amplitude self-stabilization, analogous to phase locking in Phase-Locked Loop (PLL) systems.

In addition, we provide a constructive theorem that yields an iterative scheme for approximating solutions of the corresponding boundary-value problem.

\section{\label{sec:level1} Introduction }

Nonlinear oscillators, first studied by Van der Pol~\cite{VanderPol1927}, revealed the phenomenon of frequency demultiplication. While the dynamics of single oscillators are well-understood, the behavior of interconnected networks remains an active research area. Coupled Van der Pol systems have been used to model collective synchronization in engineered and biological systems~\cite{Zontag, FitzHugh1961}. Modern trends in network science and neuromorphic computing highlight the relevance of analyzing such \emph{interconnected systems}, including recurrent neural networks (RNN) and Cohen-Grossberg networks \cite{Cohen2}, which can be formally represented as networks of nonlinear oscillators.

Oscillators \cite{Amore} and synchronization phenomena constitute a fundamental paradigm in nonlinear dynamics  \cite{Pikovsky}, \cite{Winfree}, \cite{Strog} across classical nonlinear dynamics \cite{Izhik1}, \cite{Izhik2}, \cite{Kuram}, \cite{Golub}, but the quantum regime reveals qualitatively new behaviors that defy classical intuition. The quantum analogue of the van der Pol oscillator \cite{Lee1} (with master equation for the density matrix $\rho$) — a paradigmatic self-sustained system with nonlinear damping — has emerged as a minimal model for studying quantum limit cycles, phase coherence, and quantum synchronization in driven-dissipative open systems. 

As shown in this work, the proposed method remains applicable when the density operator $\rho$ is matrix-valued (or, more generally, an operator matrix). In this setting, the commutator is defined as
$$
[H, \rho] = H\rho - \rho H,
$$
The problem can be reduced to a standard form by employing the Kronecker product representation. This transformation leads to an operator equation for a vector-valued function, where the density acts via left multiplication, and results in an effective twofold reduction of dimensionality (science).

Recent works in Applied Physics Letters have demonstrated the integration of physics-informed neural networks into the modeling of complex physical systems, where neural architectures are constrained by governing physical laws \cite{apl_pinn_adaptation_2026}.

Physics-informed neural networks have been successfully applied to nonlinear physical systems, including phase-field dynamics and coupled physical processes.

More recently, the general framework developed in semiclassical phase reduction theory for quantum synchronization extends classical phase-reduction techniques to dissipative quantum limit-cycle oscillators, showing that under certain semiclassical approximations the phase dynamics of a quantum oscillator can be accurately mapped to an effective classical stochastic equation. This bridge enables quantitative analysis of synchronization, frequency entrainment, and stability in quantum systems using tools long established in classical nonlinear dynamics. 

The analytical framework developed in this paper is not restricted to classical Van der Pol networks. Since the problem is formulated at the level of an abstract equation in a separable Hilbert space, the approach applies to a broad class of systems. In particular, depending on the specific realization, it leads to models ranging from classical to quantum generalizations of the Van der Pol equation.
After expansion in a Fourier basis, the system reduces to a nonlinear system of ordinary differential equations with a double index structure: one index runs over the natural numbers, while the other takes values in a finite set of integers, resulting in a countable system. The analysis of such systems provides a useful framework for describing amplitude and phase dynamics, including in experimental realizations of quantum Van der Pol oscillators~\cite{LIE}.

Such phenomena highlight the nontrivial interplay between coherence, dissipation, 
and nonlinearity in quantum many-body and networked oscillator systems, with 
potential applications in quantum metrology, quantum control, and neuromorphic 
quantum information processing.

We emphasize that our approach is not intended to fully describe quantum dynamics, but rather to complement existing models by providing additional analytical tools for studying amplitude constraints, phase relations, and collective behavior in weakly nonlinear regimes.

Given this emerging landscape, quantum synchronization via van der Pol-type oscillators provides a compelling framework to explore how nonlinearity, dissipation and quantum fluctuations conspire to produce robust phase coherence, nonclassical correlations and novel steady states. In the context of infinite-dimensional (Hilbert space) systems with many modes — analogously to generalized van der Pol equations in functional spaces — these results motivate the study of collective quantum oscillations, mode coupling, and the quantum-geometric structure of the corresponding phase space (we can observe some kinds of symmetry \cite{Collins1}).

We use the well known Lyapunov-Schmidt methods using of Moore-Penrose pseudoinverse operators \cite{Penr} and Lyapunov-Poincare equations ~\cite{Author2025_P}, \cite{Malkin}. 

\section{Statement of the problem}

Consider $n$ coupled Van der Pol oscillators in a separable Hilbert space $\mathcal{H}$:
\begin{equation} \label{eq:1}
\ddot{y}_m(t,\varepsilon) + Ty_m(t,\varepsilon) = \varepsilon (1 - \sum_{l = 1}^{n}\|y_l(t,\varepsilon)\|^2) \dot{y}_m(t,\varepsilon), 
\end{equation}
with boundary conditions
\begin{equation} \label{eq:2}
y_m(0,\varepsilon) - y_m(w,\varepsilon) = 0, ~~~~~~m = 1,\dots,n,
\end{equation}
\begin{equation} \label{eq:3}
\quad \dot{y}_m(0,\varepsilon) - \dot{y}_m(w,\varepsilon) = 0, \quad m=1,\dots,n.
\end{equation}

Here $T$ is an unbounded operator with bounded inverse, and $\{e_k\}_{k\in N}$ is an orthonormal basis such that
\begin{equation}
Ty_m(t,\varepsilon) = \sum_{k=1}^{+ \infty} \lambda_k c_m^k(t, \varepsilon) e_k,
\end{equation}
where $\lambda_k$ are spectral values and $c_m^k(t, \varepsilon)$  are Fourier coefficients. Unlike the one-dimensional classical Van der Pol equation, in our formulation the coefficient $T$
is an operator acting on a spatial variable, such as the Laplacian. This allows the system to support distributed or infinite-dimensional modes, which can be analyzed through corresponding amplitude equations for each mode. Expanding $y_i(t,\varepsilon)$ in this basis yields a system of nonlinear amplitude equations. 
It should be noted that we can also consider the interconnected system of partial differential equations (hyperbolic equations) in the form
$$
\frac{\partial^{2} y_{m}(t, x, \varepsilon)}{\partial t^{2}} = \Delta y_{m}(t, x, \varepsilon) + 
$$
$$
+ \varepsilon \left( 1 - \sum_{j = 1}^{n}y_{j}^{2}(t, x, \varepsilon) \right)\frac{\partial y_{m}(t, x, \varepsilon)}{\partial t}
$$
with different types of boundary conditions (as an example of an abstract model). Here $\Delta$ is the Laplacian ($\Delta y_{m}(t, x, \varepsilon) = \frac{\partial^{2}}{\partial x^{2}} y_{m}(t, x, \varepsilon))$ in the case
 when $x \in A \subset R$  and $\Delta y_{m}(t, x, \varepsilon) = \sum_{k = 1}^{l}\frac{\partial^{2}}{\partial x_{k}^{2}}y_{m}(t, x, \varepsilon)$ in the case when $x \in A \subset R^{l}$ ($A$ is a corresponding domain in $R^{l}$).

In the paper \cite{Author2025_P} (where Schr$\ddot{o}$dinger equation was also considered) the authors investigated only single oscillator of this form 
$$
\ddot{y}(t, \varepsilon) + Ty(t, \varepsilon) = \varepsilon(1 - ||y(t, \varepsilon)||^{2})\dot{y}(t, \varepsilon)
$$
with periodic conditions 
\begin{equation}
y(0,\varepsilon) - y(w,\varepsilon) = 0, 
\end{equation}
\begin{equation}
\quad \dot{y}(0,\varepsilon) - \dot{y}(w,\varepsilon) = 0.
\end{equation}

For this oscillator, the analogue of Lyapunov-Poincarre equations take the form of coupled cubic equations over a countable set of variables $(c_1^k, c_2^k)$:
\begin{equation} \label{eqs:1}
\begin{array}{ccccc}(c_1^k)^3 + 2 \sum_{j\neq k} (c_1^k (c_1^j)^2 + c_1^k (c_2^j)^2) + 3c_1^k (c_2^k)^2 - 4 c_1^k &= 0, \\
(c_2^k)^3 + 2 \sum_{j\neq k} (c_2^k (c_1^j)^2 + c_2^k (c_2^j)^2) + 3(c_1^k)^2 c_2^k - 4 c_2^k &= 0.
\end{array}
\end{equation}
Our analysis \cite{Author2025_P} showed that nonzero modes of amplitudes are restricted to a finite subset 
\begin{equation}
(c_1^k)^2 + (c_2^k)^2 = R^2,
\end{equation}
forming an $n$-dimensional torus independent of the network size. Here $c_{1}^{k}, c_{2}^{k}$	are the amplitude and phase parameters of the 
k-th nonlinear Van der Pol mode in the 
corresponding coordinate of the system, $A_{k} = \sqrt{(c_{1}^{k})^{2} + (c_{2}^{k})^{2}}$ are corresponding amplitudes. 

In this paper, it will be shown that for a system of $n$ oscillators, the equations for the generating amplitudes modes have a completely similar construction and are thus invariant with respect to the number of oscillators (i.e., for an arbitrary fixed number of oscillators that are connected to each other, the necessary condition for the existence of solutions is the condition on the amplitudes modes of oscillations, which must also lie on n-dimensional toroidal manifolds). Moreover, we can solve our system completely and explicitly identify the exact set of points on the tori corresponding to the amplitude modes, which was not done in our earlier work on this topic \cite{Author2025_P}.


We seek the solution of our system in the form
$$
y_{m}(t, \varepsilon) = \sum_{k = 1}^{+\infty}c_{m}^{k}(t, \varepsilon)e_{k}.
$$
Denote by $x_{k}^{m}(t, \varepsilon) = c_{m}^{k}(t, \varepsilon)$. After substitution in our system we can obtain the following countable systems of ordinary nonlinear equations with corresponding boundary conditions
\begin{equation} \label{Dec:140}
\left\{
\begin{aligned}
\dot{x}_{k}^{m}(t,\varepsilon) &= \sqrt{\lambda_k}\, y_k^m(t,\varepsilon), \\
\dot{y}_{k}^{m}(t,\varepsilon) &= -\sqrt{\lambda_k}\,x_k^m(t,\varepsilon)
 + \\
 + \varepsilon \left(1 - \sum_{l = 1}^{n}\sum_{j=1}^{\infty} (x_j^l(t,\varepsilon))^{2}\right)y_k^m(t,\varepsilon),
\end{aligned}
\right.
\end{equation}
with boundary conditions 
\begin{equation} \label{Dec:133}
\left\{ \begin{array}{cccccc}
x_{k}^{m}(0, \varepsilon) - x_{k}^{m}(w, \varepsilon) = 0, ~~k = 1, \dots, \infty;~~ m = 1, \dots, n, \\
y_{k}^{m}(0, \varepsilon) - y_{k}^{m}(w, \varepsilon) = 0, ~~k = 1, \dots, \infty; ~~m = 1, \dots, n.
\end{array} \right.
\end{equation}
Here $\varepsilon << 1 $ is a small parameter (so called weak coupling strength)\cite{Izhik1}, \cite{Izhik2}. 
We consider the so-called resonance case when $\lambda_{k} = \frac{4 \pi^{2}k^{2}}{w^{2}}, i \in N$ and $w = 2\pi$.
For $\varepsilon = 0$ we obtain the following generated linear boundary-value problem
\begin{equation} \label{Dec:131}
\left\{ \begin{array}{ccccc}\dot{x}_{k}^{m}(t, 0) =  k y_k^{m}(t, 0), \\
\dot{y}_{k}^{m}(t, 0) = - k x_k^{m}(t, 0),
\end{array}\right. 
\end{equation}
with boundary conditions
\begin{equation} \label{Dec:132}
\left\{ \begin{array}{cccccc} x_{k}^{m}(0, 0) - x_{k}^{m}(2\pi, 0) = 0, ~~k = 1, \dots, \infty;~~ m = 1, \dots, n, \\
y_{k}^{m}(0, 0) - y_{k}^{m}(2\pi, 0) = 0, ~~k = 1, \dots, \infty;~~ m = 1, \dots, n.
\end{array}
\right.
\end{equation}
We seek solutions $x_{k}^{m}(t, \varepsilon), y_{k}^{m}(t, \varepsilon)$ (\ref{Dec:140}), (\ref{Dec:133}) which turn (for $\varepsilon = 0$) in one of solutions of the linear generated boundary-value problem  (\ref{Dec:131}), (\ref{Dec:132}).

We seek the solution of our system in the form
$$
y_{m}(t, \varepsilon) = \sum_{k = 1}^{+\infty}c_{m}^{k}(t, \varepsilon)e_{k}.
$$
Denote by $x_{k}^{m}(t, \varepsilon) = c_{m}^{k}(t, \varepsilon)$. After substitution in our system we can obtain the following countable systems of ordinary nonlinear equations with corresponding boundary conditions
\begin{equation} \label{Dec:140}
\left\{
\begin{aligned}
\dot{x}_{k}^{m}(t,\varepsilon) &= \sqrt{\lambda_k}\, y_k^m(t,\varepsilon), \\
\dot{y}_{k}^{m}(t,\varepsilon) &= -\sqrt{\lambda_k}\,x_k^m(t,\varepsilon)
 + \\
 + \varepsilon \left(1 - \sum_{l = 1}^{n}\sum_{j=1}^{\infty} (x_j^l(t,\varepsilon))^{2}\right)y_k^m(t,\varepsilon),
\end{aligned}
\right.
\end{equation}
with boundary conditions 
\begin{equation} \label{Dec:133}
\left\{ \begin{array}{cccccc}
x_{k}^{m}(0, \varepsilon) - x_{k}^{m}(w, \varepsilon) = 0, ~~k = 1, \dots, \infty;~~ m = 1, \dots, n, \\
y_{k}^{m}(0, \varepsilon) - y_{k}^{m}(w, \varepsilon) = 0, ~~k = 1, \dots, \infty; ~~m = 1, \dots, n.
\end{array} \right.
\end{equation}
Here $\varepsilon << 1 $ is a small parameter (so called weak coupling strength)\cite{Izhik1}, \cite{Izhik2}. 
We consider the so-called resonance case when $\lambda_{k} = \frac{4 \pi^{2}k^{2}}{w^{2}}, i \in N$ and $w = 2\pi$.
For $\varepsilon = 0$ we obtain the following generated linear boundary-value problem
\begin{equation} \label{Dec:131}
\left\{ \begin{array}{ccccc}\dot{x}_{k}^{m}(t, 0) =  k y_k^{m}(t, 0), \\
\dot{y}_{k}^{m}(t, 0) = - k x_k^{m}(t, 0),
\end{array}\right. 
\end{equation}
with boundary conditions
\begin{equation} \label{Dec:132}
\left\{ \begin{array}{cccccc} x_{k}^{m}(0, 0) - x_{k}^{m}(2\pi, 0) = 0, ~~k = 1, \dots, \infty;~~ m = 1, \dots, n, \\
y_{k}^{m}(0, 0) - y_{k}^{m}(2\pi, 0) = 0, ~~k = 1, \dots, \infty;~~ m = 1, \dots, n.
\end{array}
\right.
\end{equation}
We seek solutions $x_{k}^{m}(t, \varepsilon), y_{k}^{m}(t, \varepsilon)$ (\ref{Dec:140}), (\ref{Dec:133}) which turn (for $\varepsilon = 0$) in one of solutions of the linear generated boundary-value problem  (\ref{Dec:131}), (\ref{Dec:132}). 
The set of solutions of this boundary-value problem has the form
$$
x_{k}^{m}(t, 0) = \cos(kt) c_{1 m}^{k} + \sin(kt)c_{2 m }^{k} = 
$$
$$
= \sqrt{(c_{1m}^{k})^{2} + (c_{2m}^{k})^{2}}\cos(kt - \varphi_{k}^{m}),
$$
$$
y_{k}^{m}(t, 0) = -\sin(kt) c_{1 m}^{k} + \cos(kt)c_{2 m}^{k} = 
$$
$$
= \sqrt{(c_{1m}^{k})^{2} + (c_{2m}^{k})^{2}}sin(kt - \varphi_{k}^{m}), ~~ k \in N,
$$
for any $c_{1m}^{k}, c_{2m}^{k}$.

\section{Main results}
Necessary condition on the existence of solutions of the nonlinear problem can be obtained using operator $B_{0}$ which described as in the paper \cite{Author2025_P}. After some calculations we obtain the system of Lyapunov-Poincare \cite{Malkin} equations in the form:
\begin{equation*} \label{esul:1}
\left( \begin{array}{cccccc} 
F_{i}^{1}(c_{1m}^{1}, c_{2m}^{1}, ..., )\\
F_{i}^{2}(c_{1m}^{1}, c_{2m}^{1}, ..., ) 
\end{array}\right) : = \int_{0}^{2\pi}\left(\begin{array}{ccccc} \cos(k\tau) & -\sin(k\tau) \\
\sin(k\tau) & \cos(k\tau)\end{array}\right) \times 
\end{equation*}
\begin{equation}
\times \left(\begin{array}{cccccc} 0 \\
\left(1 - \sum_{l = 1}^{n}\sum_{j = 1}^{+\infty}(x_{j}^{l}(\tau, 0))^2\right)y_{m}^{k}(\tau, 0)\end{array}\right)d\tau.
\end{equation}

Finally, we obtain 
\begin{equation} \label{eq:110}
\begin{array}{ccccc}
-4c_{1 m}^k + 2 \sum_{l = 1}^{n}\sum_{j\neq k} (c_{1 m}^k (c_{1 l}^j)^2 + c_{1 m}^k (c_{2 l}^j)^2) + \\
+ c_{1 m}^k \sum_{l = 1}^{n}\left((c_{1 l}^k)^2 + 3(c_{2l}^{k})^{2}\right) - 2\sum_{l = 1}^{n}c_{2m}^{k}c_{1l}^{k}c_{2l}^{k} &= 0, \\
-4c_{2 m}^k + 2 \sum_{l = 1}^{n}\sum_{j\neq k} (c_{2 m}^k (c_{1 l}^j)^2 + c_{2 m}^k (c_{2l}^j)^2) + \\
+  c_{2 m}^k \sum_{l = 1}^{n}\left(3(c_{1 l}^k)^2 + (c_{2l}^{k})^{2}\right) - 2\sum_{l = 1}^{n}c_{1m}^{k}c_{1l}^{k}c_{2l}^{k} &=  0.
\end{array}
\end{equation}
In the case when $m = n = 1$ we obtain result from \cite{Author2025_P}.  
From these equalities we obtain the simple necessary condition of the existence of solutions of the nonlinear boundary-value problem.

{\bf Lemma.} (amplitude constraint).  {\it 
If the periodic boundary-value problem for the coupled Van der Pol network is solvable, then only a finite number of Fourier modes are active. Moreover, for each active mode the amplitude coefficients satisfy the invariant relation
\begin{equation} \label{Dec:13}
R^{2} = \sum_{l = 1}^{n}\left((c_{1 l}^{k})^{2} + (c_{2 l}^{k})^{2}\right) = \frac{4}{2N + 1}, 
\end{equation}
where $R > 0$ is a constant independent of the mode index $k$,  $R = \frac{2}{\sqrt{2N + 1}}$. Consequently, all admissible states lie on an $nP$-dimensional toroidal manifold. (see Appendix). Moreover $c_{1m}^{k} = c_{2m}^{k}$ or $c_{1m}^{k} = -c_{2m}^{k}$.
}
 
 In this context, nonlinear oscillator networks such as those based on the Van der Pol model provide a versatile classical platform for exploring the emergence of phase quantization, mode hybridization, and synchronization. In the weakly nonlinear regime considered here, we show that the generating mode amplitudes satisfy a discrete phase condition 
$c_{1m}^{k} = c_{2m}^{k}$ or $c_{1m}^{k} = -c_{2m}^{k}$, which can be regarded as a spinor-like degree of freedom, analogous to two-component band states in quantum materials. Thus, our model naturally introduces a geometrically constrained state space where different branches are separated by a phase inversion — reminiscent of topologically distinct quantum phases.
This connection suggests that the Van der Pol lattice can serve as a classical analogue of systems exhibiting quantum-geometric effects. Upon quantization of collective modes (phonons, polaritons, superconducting resonators, etc.), the same structure may lead to the realization of Berry-phase-induced response, topologically protected synchronization, or geometric control of collective oscillations. Consequently, the framework developed here offers both a new perspective on nonlinear dynamics and a promising foundation for exploring topological collective behavior in engineered platforms bridging classical and quantum condensed matter physics. Finally, we have the following assertion.

{\bf Theorem 1.} (necessary condition). {\it 
In addition to the amplitude constraint stated in Lemma 1, the Lyapunov–Poincaré reduction yields an additional phase condition for all nonzero modes:
$$
c_{1m}^{k} = \pm c_{2m}^{k}.
$$
Equivalently, in polar coordinates 
$$
c_{1m}^{k} = R \cos\varphi_{m}^{k}, c_{2m}^{k} = R \sin \varphi_{m}^{k},
$$ 
the admissible phases are restricted to
$$
\varphi_{m}^{k} \in \left\{ \frac{\pi}{4}, \frac{3\pi}{4}, \frac{5\pi}{4}, \frac{7\pi}{4} \right\}
$$
Therefore, the continuous toroidal manifold is reduced to a finite discrete set of configurations:
$$
\left(c_{1m}^{k}, c_{2m}^{k} \right) \in \left\{ \left( \pm \frac{R}{\sqrt{2}}, \pm \frac{R}{\sqrt{2}}\right), \left( \pm \frac{R}{\sqrt{2}}, \mp \frac{R}{\sqrt{2}}\right)\right\},
$$
$m = 1, \dots, n$.
This set consists of $2N + 1$ discrete phase-locked configurations lying on the invariant toroidal manifold.
}

This result describes a symmetry-breaking mechanism on the invariant torus, leading to phase-locked discrete states analogous to spinor-type degeneracy in coupled oscillator systems.

{\bf Remark 1.} {\it As you can see that the form of equations for generating parameters (\ref{eqs:1}) has the same form as (\ref{eq:110}) and accordingly is an invariant for interconnected systems with any number of oscillators.}

{\bf Remark 3.} {\it It should be emphasized that, if the number of parameters $c_{1m}^{k}, c_{2m}^{k}$ is initially known to be finite, then the necessary and sufficient conditions coincide. In this case, we obtain a criterion for the existence of a unique periodic solution to the coupled Van der Pol system of nonlinear ordinary differential equations in the form
\begin{equation}
\left\{
\begin{aligned}
\dot{x}_{k}^{m}(t,\varepsilon) &= k\, y_k^m(t,\varepsilon), \\
\dot{y}_{k}^{m}(t,\varepsilon) &= -k\,x_k^m(t,\varepsilon)
 + \\
 + \varepsilon \left(1 - \sum_{l = 1}^{n}\sum_{j=1}^{P} (x_j^l(t,\varepsilon))^{2}\right)y_k^m(t,\varepsilon)
\end{aligned}
\right.
\end{equation}
with periodic conditions.
}

{\bf Remark 4.} {\it Moreover, the final theorem also works in the case when the number of oscillators are infinite ($n = +\infty$) (invariant has the same form). See also \cite{Phon}, \cite{Liu}, \cite{Timokha}.}

In conclusion, we discuss the amplitude–frequency characteristics and their connection to the analytical results. Now, we can represent sufficient condition of the existence of solutions of the periodic boundary-value problem with algorithm of finding of approximative solution.

{\bf Theorem 2.} (sufficient condition). {\it Under condition of the previous theorem periodic problem is solvable. The sequence of the corresponding solutions can be found using of the following iterative processes
$$
\left\{ \begin{array}{ccccc}
\overline{Z}_{l + 1}^{m}(t, \varepsilon) = \varepsilon F(Z_{l}^{m}(t, 0) + H_{l}^{m}(t, \varepsilon)), \\
Z_{l + 1}^{m}(t, \varepsilon) = U^{m}(t)P_{N(B_{0})}\overline{c}_{m} + \varepsilon \overline{Z}_{l + 1}^{m}(t, \varepsilon),
\end{array} \right.
$$
where $Z_{l + 1}^{m}(t, \varepsilon) = Z_{l}^{m}(t, \varepsilon) + H_{l}^{m}(t, \varepsilon)$ and 
$$
F(Z_{l}^{m}(t, 0) + H_{l}^{m}(t, \varepsilon)) =
$$
$$
= U^{m}(t)B_{0}^{+}\int_{0}^{T}U^{m}(T)U^{-1 m}(\tau)G_{m}^{l}(\tau, \varepsilon)d\tau + 
$$
$$
+ \int_{0}^{t}U^{m}(t)U^{-1 m}(\tau)G_{m}^{l}(\tau, \varepsilon)d\tau,
$$
where 
$$
G_{m}(t, \varepsilon) = \left( \begin{array}{ccccc} 0  \\
(1 - \sum_{l = 1}^{n}\sum_{j = 1}^{\infty}(x_{j}^{l}(t, \varepsilon))^{2}y_{1}^{m}(t, \varepsilon)) \\
0 \\
(1 - \sum_{l = 1}^{n}\sum_{j = 1}^{\infty} (x_{j}^{l}(t, \varepsilon))^{2}y_{2}^{m}(t, \varepsilon)) \\
0 \\
...
\end{array} 
\right),
$$
$$
U^{m}(t) = \left(\begin{array}{ccccc} 
J_{2}(\lambda_{1}t) & 0 &\cdots & 0 & \cdots\\
0 & J_{2}(\lambda_{2}t) & \cdots & 0 & \cdots \\
\cdots & \cdots & \cdots & \cdots & \cdots \\
0 & 0 & 0 & 0 & J_{2}(\lambda_{n}t)  \\
\dots & \dots & \dots & \dots & \dots 
\end{array} \right),
$$
$$
J(\lambda t) = \left( \begin{array}{cccc} 
cos(\lambda t) & \sin(\lambda t) \\
-sin(\lambda t) & cos(\lambda t)  
\end{array} \right),
$$
$$
\lambda_{k} = k, w = 2\pi,
$$
$$
B_{0} = I - U_{m}(T), ~~~P_{N(B_{0}^{*})}U^{m}(T) = 0.
$$
}
The measured amplitude-frequency behavior corresponds to classical PLL-like characteristics of a Van der Pol generator, illustrating how nonlinear resistance stabilizes the amplitude: gain exceeds losses at small amplitudes, while losses dominate at large amplitudes, yielding a stable limit cycle.
In our analysis of coupled Van der Pol networks, we demonstrate the existence of an invariant amplitude structure:
$$
(c_{1m}^{k})^{2} + (c_{2m}^{k})^{2} = r^{2},
$$
which in a classical triode circuit corresponds to a stable oscillation amplitude observed at the peak frequency. For a single oscillator, the condition $c_{1}^{2} + c_{2}^{2} = 1$ reproduces classical amplitude stability.
From an engineering perspective, $c_{1}$ and $c_{2}$ correspond to harmonic components of voltage and current, or phase-space coordinates of the oscillator. The toroidal sets obtained for networks indicate that even in complex systems there exists a restricted set of stable oscillatory modes, analogous to the amplitude stabilization of a single triode.
Experimental verification remains an open task. In a laboratory network of multiple Van der Pol circuits (e.g., RC oscillators with tubes or transistors), the invariant structure $(c_{1m}^{k})^{2} + (c_{2m}^{k})^{2} = r^{2}$ should correlate with stable amplitudes, and resonance peaks and phase shifts can be compared with the theoretical values $c_{1m}^{k}, c_{2m}^{k}$.

An essential feature of the considered nonlinear system is that each oscillator is coupled to every other one, forming a fully connected network (analogous to a complete graph). This all-to-all interaction imposes strong constraints on the collective dynamics, so that the amplitudes of all modes must satisfy a shared normalization condition
(\ref{Dec:13}) which ensures that any change in one oscillator immediately affects the allowed amplitudes of all others. As a result, the network admits only a finite set of discrete, geometrically distinct limit-cycle configurations on the toroidal manifold, rather than a continuous distribution.

The proposed framework is consistent with related approaches that have been successfully applied to complex systems in recent high-impact studies \cite{Humanity} 
(including the present author), demonstrating its applicability to problems of comparable complexity.

The analytical invariant structures identified in this work are not restricted to traditional physical applications. 
Recent studies in physics-informed machine learning demonstrate that data-driven methods can discover conserved quantities and symmetry invariants from trajectory data in nonlinear dynamical systems, highlighting the role of constraint-based representations across different domains~\cite{Ha2021}. 

Beyond physical applications, these ideas may also be relevant for structured high-dimensional problems arising in modern machine learning, where invariant-based representations play a key role in data-driven modeling. The proposed iterative scheme is complemented by numerical simulations, 
illustrating the trajectories and phase portraits of the coupled system (see Fig.~1). 
These results demonstrate the formation of stable oscillatory regimes and provide 
computational evidence for the invariant amplitude structure derived analytically. 
In particular, the dynamics remain confined to toroidal manifolds, confirming the 
predicted amplitude self-stabilization mechanism characteristic of Van der Pol–type systems.

Trajectories you can find from the algorithm (example).

\begin{figure}[h!]
    \centering
    \includegraphics[width=1\textwidth]{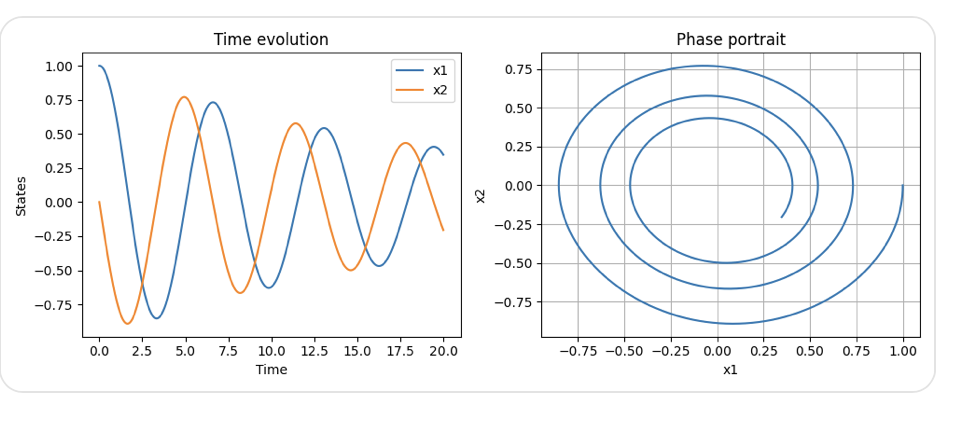} 
    \caption{Trajectories and phase portrait}
    \label{fig:my_label}
\end{figure}

It should be noted that the main result doesn't depend from the number $n$ ($n$ can be also infinite). In such a way we can investigate not only classical van der Pol model  but different kinds of quantuum model.

\section{An Applied Example of Coupled Van der Pol Networks: Line-Level OCR}

We take the coupled Van der Pol network framework of Theorem~1 and turn it into a differentiable layer inside a deliberately simple deep OCR model. This applied example has two aims:
\begin{enumerate}
    \item show that the framework can be embedded as a plug-in NN layer with negligible parameter overhead and standard end-to-end CTC training;
    \item show that its amplitude-stabilizing dynamics behave as a feature restorer for corrupted inputs, recovering recognition accuracy on document-level degradations (blur, low contrast, historical wear, ink bleed, noise) where the same CRNN without the layer degrades sharply.
\end{enumerate}

The task is line-level handwriting recognition on the \texttt{Teklia/IAM-line}~\cite{IAM_DATASET} dataset from HuggingFace (grayscale handwriting lines, height-normalized to 128~px), and the baseline is a raw standard CRNN with three 3×3 conv blocks $(64{\to}128{\to}256$, MaxPool x2 after each), height collapsed by adaptive average pooling, a two-layer bidirectional LSTM (hidden 256), and a linear feed-froward layer + CTC head. 

{\bf Architecture}. We insert a differentiable Van der Pol (VDP) stabilization layer inside this baseline after \texttt{conv2} (feature map $B\times128\times32\times W{/}4$). Features are projected to $K=32$ quadrature coordinates $(x_k,y_k)$ by two 1×1 convolutions and L2-normalized onto a per-channel learnable radius $R_\theta > 0$ ($N=8$ groups, $g=4$). Two distinct constants govern the layer: (i) Theorem 1's fixed theoretical group-energy target $R^2_\text{paper}=4/(2N+1)\approx 0.235$, which appears only in the amplitude regularizer $L_\text{amp}=\bigl(\sum_{k\in G}(x_k^2+y_k^2)-R^2_\text{paper}\bigr)^2$ pinning the group-summed energy during training; and (ii) the learnable per-channel radius $R_\theta$, initialized to $R_\text{init}=\sqrt{R^2_\text{paper}/g}=0.2425$ and used only in the projection normalization — writing $(\tilde x_k,\tilde y_k)=(W_x f, W_y f)_k$ for the raw $1{\times}1$-conv output,
$$
(x_k,y_k)\;=\;R_\theta\cdot\frac{(\tilde x_k,\tilde y_k)}{\|(\tilde x_k,\tilde y_k)\|_2},
$$
so every $(x_k,y_k)$ pair sits exactly on the per-channel radius-$R_\theta$ sphere at the start of the ODE integration (the ODE itself never references $R_\theta$). Five Heun (RK2) steps then integrate
$$
\dot x_k=\omega_k y_k,\quad \dot y_k=-\omega_k x_k+\varepsilon_k(1-E_k)y_k+(\mathcal L x)_k,
$$

Each mode $k$ is a harmonic oscillator with learnable angular frequency
$\omega_k>0$: the linear pair $(\omega_k y_k,\,-\omega_k x_k)$ alone
conserves $x_k^2 + y_k^2$, so $(x_k,y_k)$ moves round a circle at speed
$\omega_k$. The Van der Pol term $\varepsilon_k(1-E_k)\,y_k$ turns that
passive circle into an attractor: when the local energy
$E_k = \sum_{k'\in G}\bigl(\kappa \ast x_{k'}^2\bigr)$ is below the
invariant, energy is injected; when it exceeds the invariant, energy is
dissipated, so the trajectory is pulled back to the $R$-sphere at every
step. Here $\kappa$ is a fixed, non-trainable $1\times21$ box kernel ---
a same-size mean along the width axis --- that replaces the paper's
global energy sum by a local neighborhood mean at the cost of one
depthwise convolution. 

The coupling term
$(\mathcal L x)_k = \sum_q A_{pq}(x_q - x_p)$ is a learnable
depthwise, same-padded $1\times21$ convolution along the width axis ---
one kernel of length $21$ per channel. Because each IAM character spans roughly ten feature columns at
\texttt{conv2} resolution ($W/4$), the $21$-column window couples a
character to its immediate left/right neighbors --- "adjacent-site"
pairing of Theorem~1. 

Together,
$(\omega_k, \varepsilon_k, \mathcal L)$ realize the coupled
Van der Pol network on the CNN feature grid. VDP stabilization layer's
output is combined with the projection through a convex residual
$x = (1{-}\alpha)\,x_{\text{pre}} + \alpha\,x_{\text{post}}$ with
$\alpha = \sigma(\tilde\alpha) \in (0,1)$
($\sigma$ is the sigmoid).

{\bf Training}. Both models are trained end-to-end for 50 epochs on the IAM-line training split ($6{,}482$ lines, batch size $24$) with the Adam optimizer
(learning rate $10^{-3}$, cosine annealing to zero) on two NVIDIA~T4 GPUs. VDP layer optimizes the joint loss
$L = L_{\text{CTC}} + \lambda_{\text{amp}}\, L_{\text{amp}}$ with
$\lambda_{\text{amp}} = 0.01$; the phase regulariser $L_{\text{phase}}$
associated with Theorem~1's $c^1_{km} = \pm c^2_{km}$ condition is left
disabled ($\lambda_{\text{phase}} = 0$), as amplitude alone was
sufficient for the invariant to hold.

{\bf Results}. Robustness on the test set is assessed on nine document-level corruptions: blur, Gaussian noise, JPEG, low contrast, rotation, elastic warp, broken strokes, ink bleed, and historical degradation.

\begin{table}[h]
\caption{\label{tab:headline}Character (CER) and word (WER) error rates on the IAM-line test set (2915 lines). $\Delta$ is the change from clean to mean corrupted.}
\centering
\setlength{\tabcolsep}{10pt}
\begin{tabular}{lccccccc}
\toprule
Model & Params
    & \multicolumn{2}{c}{Clean$\downarrow$}
    & \multicolumn{2}{c}{Mean corrupted$\downarrow$}
    & \multicolumn{2}{c}{Degradation $\Delta \downarrow$} \\
 & & CER & WER & CER & WER & CER & WER \\
Baseline     & 3.55\,M          & 0.1822          & 0.4627          & 0.2933          & 0.6336          & $+0.1111$          & $+0.1709$          \\
\textbf{VDP} & +9.3\,K & \textbf{0.1724} & \textbf{0.4475} & \textbf{0.2629} & \textbf{0.5839} & $\mathbf{+0.0904}$ & $\mathbf{+0.1364}$ \\
\end{tabular}
\end{table}

The mean corrupted CER drops by $10.4\%$ and the CER degradation gap between clean and corrupted inputs shrinks by $19\%$ ($0.111 \to 0.090$). Table 2 breaks this down across the nine document-level corruptions and matches the theoretical prediction: the largest restorations occur on smooth, amplitude-attenuating corruptions — blur ($-25.2\%$), historical wear ($-20.0\%$), low contrast ($-19.8\%$) --- exactly the class the amplitude invariant $\sum_{k\in G}(x_k^2+y_k^2)=R^2_\text{paper}$ is designed to counteract, since these attenuate the projected feature magnitudes without destroying structure. Gaussian noise ($-10.6\%$) and ink bleed ($-8.9\%$) show smaller but consistent restoration. Only elastic warp regresses ($+5.1\%$), consistent with the width-local coupling assumption breaking under non-affine spatial distortion.

At the end of training the learnable radius converges to $R_\theta = 0.2303$, within $5\%$ of the theoretical initialization $R_\text{init} = 0.2425$, and the convex mixing weight to $\alpha = 0.173$. In other words the network keeps the amplitude invariant close to the value predicted by Theorem 1 and actively applies a $\sim 17\%$ VDP correction on top of the identity path (rather than collapsing $\alpha$ to $0$), confirming that the theoretical construction is used as designed and not bypassed during training.

\begin{table}[h]
\caption{\label{tab:headline}Character (CER) error rates breakdown across corrupted groups on the IAM-line test set (2915 lines).}
\centering
\begin{tabular}{lrrr}
\toprule
Corruption      & Baseline$\downarrow$ & VDP$\downarrow$  & $\Delta\%$ \\
blur            & 0.3687   & 0.2758  & $\mathbf{-25.2}$ \\
historical      & 0.2826   & 0.2260  & $\mathbf{-20.0}$ \\
low\_contrast   & 0.2814   & 0.2258  & $\mathbf{-19.8}$ \\
gaussian\_noise & 0.2909   & 0.2601  & $-10.6$ \\
ink\_bleed      & 0.3095   & 0.2818  & $-8.9$ \\
broken\_strokes & 0.2002   & 0.1879  & $-6.1$ \\
clean           & 0.1822   & 0.1724  & $-5.4$ \\
rotation        & 0.2519   & 0.2390  & $-5.1$ \\
jpeg            & 0.1848   & 0.1755  & $-5.0$ \\
elastic         & 0.4700   & 0.4939  & $\mathbf{+5.1}$ \\
\end{tabular}
\end{table}

Because this is a minimal-tuning case study on a raw baseline, the absolute numbers should be read as evidence that the Van der Pol framework can inject useful inductive bias into a deep model, not as a competitive OCR result. Comparison against tuned OCR baselines, multi-seed variance analysis, and per-component ablations of the VDP layer are left for follow-up work; the present numbers are meant only to establish that the framework is usable inside a deep architecture and behaves as the theory predicts on the corruption types that match its assumptions.

\section{APPENDIX} 
Conditions of the solvability. We can rewrite conditions of the solvability in the form of the following countable system of equations in the form:
$$
\int_{0}^{2\pi} sin (k\tau)\left(1 - \sum_{l = 1}^{n}\sum_{j = 1}^{+\infty}(x_{j}^{l}(\tau, 0))^{2} \right)y_{m}^{k}(\tau, 0)d\tau = 0,
$$
$$
\int_{0}^{2\pi} cos (k\tau)\left(1 - \sum_{l = 1}^{n}\sum_{j = 1}^{+\infty}(x_{j}^{l}(\tau, 0))^{2} \right)y_{m}^{k}(\tau, 0)d\tau = 0.
$$
After substitution we obtain the following equations
$$
\int_{0}^{2\pi}sin(k\tau)\left( - sin(k\tau)c_{1m}^{k} + cos(k\tau) c_{2m}^{k}\right) d\tau =
$$
$$
= \sum_{l = 1}^{n}\sum_{j = 1}^{+\infty}\int_{0}^{2\pi}sin(k\tau) \left(cos(j\tau)c_{1l}^{j} + sin(j\tau)c_{2l}^{j}\right)^{2} \times
$$
$$
\times \left(-sin(k\tau)c_{1m}^{k} + cos(k\tau)c_{2m}^{k}\right)d\tau,
$$
$$
\int_{0}^{2\pi}cos(k\tau)\left( - sin(k\tau)c_{1m}^{k} + cos(k\tau) c_{2m}^{k}\right) d\tau =
$$
$$
= \sum_{l = 1}^{n}\sum_{j = 1}^{+\infty}\int_{0}^{2\pi}cos(k\tau) \left(cos(j\tau)c_{1l}^{j} + sin(j\tau)c_{2l}^{j}\right)^{2} \times
$$
$$
\times \left(-sin(k\tau)c_{1m}^{k} + cos(k\tau)c_{2m}^{k}\right)d\tau
$$
or in the extended form
$$
\underline{J}_{1}^{k} c_{1m}^{k} + \underline{J}_{2}^{k}c_{2m}^{k} = \sum_{l = 1}^{n}\sum_{j = 1}^{+\infty} \underline{J}_{3}^{k, j}c_{1m}^{k}(c_{1l}^{j})^{2} + \sum_{l = 1}^{n}\sum_{j = 1}^{+\infty}\underline{J}_{4}^{k, j}(c_{1l}^{j})^{2}c_{2m}^{k} + 
$$
$$
+ \sum_{l = 1}^{n}\sum_{j = 1}^{+\infty}\underline{J}_{5}^{k, j}(c_{2l}^{j})^{2}c_{1m}^{k} + \sum_{l = 1}^{n}\sum_{j = 1}^{+\infty}\underline{J}_{6}^{k, j}(c_{2l}^{j})^{2}c_{2m}^{k}  + 
$$
$$
 + \sum_{l = 1}^{n}\sum_{j = 1}^{+\infty} \underline{J}_{7}^{k, j}c_{1l}^{j}c_{2l}^{j})c_{1m}^{k} + \sum_{l = 1}^{n}\sum_{j = 1}^{+\infty} \underline{J}_{8}^{k, j}c_{2m}^{k}c_{1l}^{j}c_{2l}^{j}.
$$

$$
\overline{J}_{1}^{k} c_{1m}^{k} + \overline{J}_{2}^{k}c_{2m}^{k} = \sum_{l = 1}^{n}\sum_{j = 1}^{+\infty} \overline{J}_{3}^{k, j}c_{1m}^{k}(c_{1l}^{j})^{2} + \sum_{l = 1}^{n}\sum_{j = 1}^{+\infty}\overline{J}_{4}^{k, j}(c_{1l}^{j})^{2}c_{2k}^{m} + 
$$
$$
+ \sum_{l = 1}^{n}\sum_{j = 1}^{+\infty}\overline{J}_{5}^{k, j}c_{1l}^{j}c_{2l}^{j} + \sum_{l = 1}^{n}\sum_{j = 1}^{+\infty}\overline{J}_{6}^{k, j}c_{1l}^{j}c_{2l}^{j}c_{2m}^{k}  + 
$$
$$
 + \sum_{l = 1}^{n}\sum_{j = 1}^{+\infty} \overline{J}_{7}^{k, j}c_{1m}^{k}(c_{2l}^{j})^{2} + \sum_{l = 1}^{n}\sum_{j = 1}^{+\infty} \overline{J}_{8}^{k,j}c_{2m}^{k}(c_{2l}^{j})^{2}, 
$$
where 
$$
\underline{J}_{1}^{k} = - \int_{0}^{2\pi}sin^{2}(k\tau)d\tau = \\
- \frac{\pi}{2},
$$
$$
\overline{J}_{1}^{k} = -\frac{1}{2}\int_{0}^{2\pi}sin(2k\tau)cos^{2}(j\tau)d\tau = -\frac{\pi}{2},   
$$
$$
\underline{J}_{2}^{k} = \frac{1}{2}\int_{0}^{2\pi} sin(2k\tau)d\tau, ~~~\overline{J}_{2}^{k} = \int_{0}^{2\pi}cos^{2}(k\tau)d\tau = \frac{\pi}{2},
$$
$$
\underline{J}_{3}^{k, j} = - \frac{1}{4}\int_{0}^{2\pi} (1 - cos(2k\tau)))(1 + cos(2j\tau)d\tau =  \\
\frac{\pi}{2},
$$
$$
\overline{J}_{3}^{k, j} = -\frac{1}{2}\int_{0}^{2\pi}sin(2k\tau)d\tau = 0, 
$$
$$
\underline{J}_{4}^{k, j} = \frac{1}{2}\int_{0}^{2\pi}sin(2k\tau)cos^{2}(j\tau)d\tau = 0,
$$
$$
\overline{J}_{4}^{k, j} = \int_{0}^{2\pi}cos^{2}(k\tau)cos^{2}(j\tau)d\tau = \frac{\pi}{2},
$$
$$
\underline{J}_{5}^{k, j} = - \int_{0}^{2\pi} sin^{2}(k\tau)sin^{2}(j\tau)d\tau = - \frac{\pi}{2},
$$
$$
\overline{J}_{5}^{k, j} = -\frac{1}{2}\int_{0}^{2\pi}sin(2k\tau)sin(2j\tau)d\tau, 
$$
$$
\underline{J}_{6}^{k, j} =  \frac{1}{2}\int_{0}^{2\pi}sin(2k\tau)sin^{2}(j\tau)d\tau,
$$
$$
\overline{J}_{6}^{k, j} = \int_{0}^{2\pi}sin(2j\tau)cos^{2}(k\tau)d\tau, 
$$
$$
\underline{J}_{7}^{k, j} = -\int_{0}^{2\pi} sin^2(k\tau)sin(2j\tau)d\tau, 
$$
$$
\overline{J}_{7}^{k, j} = -\frac{1}{2}\int_{0}^{2\pi}sin(2k\tau)sin{2}(j\tau)d\tau,
$$
$$
\underline{J}_{8}^{k, j} = \frac{1}{2}\int_{0}^{2\pi}sin(2k\tau)sin(j\tau)d\tau,
$$
$$
\overline{J}_{8}^{k, j} = \int_{0}^{2\pi}cos^{2}(k\tau)sin^{2}(j\tau)d\tau.
$$

After simplest checking we obtain (all remaining integrals vanish due to orthogonality):
$$
\overline{J}_{1}^{k} = - \int_{0}^{2\pi} \sin^2(k\tau)\, d\tau = - \pi,
\qquad
\int_{0}^{2\pi} \cos^2(k\tau)\, d\tau = \pi,
$$
$$
\int_{0}^{2\pi} \sin(k\tau)\sin(j\tau)\, d\tau =
\begin{cases}
\pi, & k = j, \\
0, & k \ne j,
\end{cases}
$$
$$
\int_{0}^{2\pi} \cos(k\tau)\cos(j\tau)\, d\tau =
\begin{cases}
\pi, & k = j, \\
0, & k \ne j,
\end{cases}
$$
$$
\overline{J}_{3}^{k, j} = - \int_{0}^{2\pi} \sin^2(k\tau)\cos^2(j\tau) = \left\{ \begin{array}{cccc} -\frac{\pi}{2}, \quad k \ne j, \\
-\frac{3\pi}{4}, \quad k = j.
\end{array} \right.
$$
$$
\overline{J}_{5}^{k, j} = -\int_{0}^{2\pi}\sin^{2}(k\tau)\sin^{2}(j\tau)d\tau = \left\{ \begin{array}{cccc} -\frac{\pi}{2}, \quad k \ne j, \\
-\frac{3\pi}{4}, \quad k = j.
\end{array} \right.
$$
Finally, we obtain the following system of equations

$$
-4c_{1m}^{k} + 2\sum_{l = 1}^{n}\sum_{i \neq j}c_{1m}^{k} \left((c_{1l}^{j})^{2} + (c_{2l}^{j})^{2}\right) +
$$
$$
+ \sum_{l = 1}^{n}(c_{1m}^{k}(c_{1l}^{k})^{2} + 3c_{1m}^{k}(c_{2l}^{k})^{2})  = 2\sum_{l = 1}^{n}c_{2m}^{k}c_{1l}^{k}c_{2l}^{k},
$$

$$
-4c_{2m}^{k} + 2\sum_{l = 1}^{n}\sum_{i \neq j}c_{2m}^{k} \left((c_{1l}^{j})^{2} + (c_{2l}^{j})^{2}\right) +
$$
$$
+ \sum_{l = 1}^{n}(3c_{2m}^{k}(c_{1l}^{k})^{2} + c_{2m}^{k}(c_{2l}^{k})^{2})  = 2\sum_{l = 1}^{n}c_{1m}^{k}c_{1l}^{k}c_{2l}^{k},
$$

If we denote $A_{j}^{2} = \sum_{l = 1}^{n}\left((c_{1l}^{j})^{2} + (c_{2l}^{j})^{2}\right) $ then
as in ... we can obtain that $A_{j}^{2} = A_{k}^{2} = A^{2}$ or $A_{k} = 0$. Moreover, the series $\sum_{j}A_{j}^{2}$ is converge if and only if the finite number of $A_{j}^{2}$ are nonzero. Let 
$$
x = \sum_{l = 1}^{n} (c_{1l}^{k})^{2}, ~~y = \sum_{l = 1}^{n}(c_{2l}^{k})^{2}, ~~x + y = A^{2}, B = \sum_{l = 1}^{n}c_{1l}^{k}c_{2l}^{k}.
$$

Then we obtain the following system
$$
c_{1m}^{k}\left(-4 + (2N + 1)A^{2}  + 2y\right) = 2c_{2m}^{k}B,
$$
$$
c_{2m}^{k}\left(-4 + (2N + 1)A^{2} +2x\right) = 2c_{1m}^{k}B.
$$
After multiplying we obtain the following system of equations
$$
c_{1m}^{k}c_{2m}^{k}(-4 + (2N + 1)A^{2}  + 2y)(-4 + (2N + 1)A^{2}  +2x) =
$$
$$
= 4c_{1m}^{k}c_{2m}^{k}B^{2}.
$$
If $c_{1m}^{k} \neq 0, c_{2m}^{k} \neq 0$, then 
\begin{equation} \label{May:20}
4B^{2} = (-4 + (2N + 1)A^{2}  + 2y)(-4 + (2N + 1)A^{2} + 2x).
\end{equation}
Multiply the first equation by $c_{1m}^{k}$ and sum over $m$ from 1 to $n$ we obtain
$$
\sum_{m = 1}^{n}(c_{1m}^{k})^{2}(-4 + (2N + 1)A^{2} + 2y) =  
$$
$$
= x(-4 + (2N + 1)A^{2} + 2y)= 2 \sum_{m = 1}^{n}c_{1m}^{k}c_{2m}^{k}B = 2B^{2}.
$$
Substituting $B^{2}$ from (17) we obtain the following equation
$$
 2x(-4 + (2N + 1)A^{2} + 2y)= (-4 + (2N + 1)A^{2}  + 2y) \times
$$
$$
\times (-4 + (2N + 1)A^{2} + 2x).
$$
Consider the case when 
$$
 -4 + (2N + 1)A^{2} + 2y \neq 0,
$$
then 
$$
2x = -4 + (2N + 1)A^{2} + 2x. 
$$
Finally, we have that 
$$
A^{2} = \frac{4}{2N + 1}
$$
and our initial system take the following form
$$
c_{1m}^{k}y = c_{2m}^{k}B,
$$
$$
c_{2m}^{k}x = c_{1m}^{k}B.
$$
From these equations we obtain that
$$
\frac{(c_{1m}^{k})^{2}}{(c_{2m}^{k})^{2}} = \frac{x}{y}.
$$
Since $B^{2} = xy$ or in the extended form
$$
\left(\sum_{l = 1}^{n}c_{1l}^{k}c_{2l}^{k}\right)^{2} = \sum_{l = 1}^{n}(c_{1l}^{k})^{2}\sum_{l = 1}^{n}(c_{2l}^{k})^{2} 
$$
we have equality in the Kochi-Bunyakovskiy inequality. It means that $\frac{c_{1p}^{k}}{c_{2p}^{k}} = \lambda, p = \overline{1, n}$ and $\lambda^{2} = \frac{x}{y}$.  Then 
$$
\left(\sum_{l = 1}^{n}c_{1l}^{k}c_{2l}^{k}\right)^{2} = \frac{1}{\lambda^{2}}\left(\sum_{l = 1}^{n}(c_{1l}^{k})^{2}\right)^{2} = \frac{x^{2}}{\lambda^{2}} = xy. 
$$
Finally
$$
y^{2} = xy.
$$
If $y \neq 0$ then $y = x$ and 
$$
\frac{(c_{1m}^{k})^{2}}{(c_{2m}^{k})^{2}} = 1.
$$
Finally 
$$
c_{1m}^{k} = \pm c_{2m}^{k}.
$$


\section{Acknowledgements}
We wish to acknowledge the support of the National Research Foundation of Ukraine 
(Project number 2025.07/0014. Project name: Modern Problems of Mathematical Analysis and Geometric Function Theory).

\end{document}